%% file: simuw_2009.tex
\documentclass[12pt]{amsart}
\usepackage{graphicx, verbatim}

\input{macros}

\begin{document}

\title{A Young Person's Guide to the Hopf Fibration}
\author{Zachary Treisman}

\maketitle

\begin{figure}
\centerline{\resizebox{\textwidth}{!}{\includegraphics{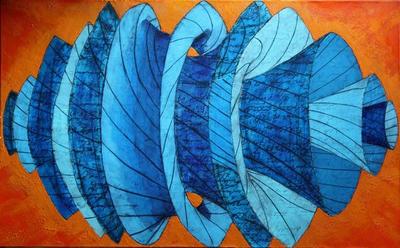}}}
\caption{ {\em Hopf Fibration} \newline Lun-Yi Tsai}
\end{figure}

The purpose of these notes is to introduce a mathematical structure which goes by the name {\em the Hopf fibration}, and demonstrates a number of surprising and beautiful things.  The Hopf fibration is a map showing a connection between two spheres, a two dimensional sphere, and a three dimensional sphere.  In order to understand this map, we will need to develop a few tools.  There will be a few detours along the way, in order to develop enough familiarity with the tools being used so that the student is sufficiently impressed by the structures that are eventually uncovered.

The Hopf fibration and the mathematics that are developed along the way makes for some very interesting visual images.  The paintings that have been included here are all the work of Lun-Yi Tsai, an artist, a mathematician, and a good friend.  The computer generated images I have produced using the programs {\em surf} and {\em jenn3d}, both of which are freely available online, as well as {\em Mathematica}.

\input{complex_numbers}

\input{spheres}

\input{quaternions}

\input{hopf_fibration}

\nocite{*}
\bibliographystyle{hamsplain}
\bibliography{simuw_2009}

\end{document}

%% file: macros.tex
\newcommand{\Pp}{{\bf P}}

\newcommand{\C}{{\mathbb C}}

\newcommand{\Q}{{\mathbb Q}}

\newcommand{\R}{{\mathbb R}}

\DeclareFontFamily{OMS}{rsfs}{\skewchar\font'60}
\DeclareFontShape{OMS}{rsfs}{m}{n}{<-5>rsfs5 <5-7>rsfs7 <7->rsfs10 }{}
\DeclareSymbolFont{rsfs}{OMS}{rsfs}{m}{n}
\DeclareSymbolFontAlphabet{\scr}{rsfs}

\newcommand{\fee}{\varphi}

\newcommand{\twobytwo}[4]{\left(\begin{array}{cc}
#1&#2\\ #3&#4
\end{array}\right)}

\theoremstyle{definition}
\newtheorem{defn}{Definition}[section]
\theoremstyle{plain}
\newtheorem{thm}[defn]{Theorem}

\theoremstyle{definition}
\newtheorem{ex}[defn]{Example}
\theoremstyle{definition}

\theoremstyle{remark}

\theoremstyle{definition}
\newtheorem{exercise}[defn]{Exercise}

%% file: complex_numbers.tex
\section{Complex Numbers}

Our first task is to introduce you to the {\em complex numbers}.  The introduction is geometric, and appeals to the intuitive basis for the real numbers as measurements that scale.  One thing that this development is not is historical.  This is because the first ways of thinking of something are not always the easiest to understand.  In this section you are often asked to show algebraic facts by drawing a picture.  Not only is this rather unusual, it is rather subtle, but it can be done with all the rigor of a traditional development.    The main idea is to see intuitively what these basic properties look like, as pictures, and to become comfortable with complex numbers.

\subsection{Fields}

The real numbers ($\R$ for short) are familarly represented on a number line, marked by important examples of real numbers, such as zero, one, and so on.  To add two numbers on the number line, we need to know where zero is.  Then we define the sum $a+b$ as the point on the line that comes from stacking the two lengths next to each other.  
\begin{exercise}Draw a picture of $2+3=5$.\end{exercise}
To multiply two numbers graphically, we need to know where zero is, and where one is as well.  Then we can define the product $ab$ as the number that $b$ reaches when the line is scaled so that 1 is at $a$. 
\begin{exercise}Draw a picture of $2\times3=6$.\end{exercise}

Frequently, a number system is abstractly defined using a set of axioms, or rules.  An undergraduate analysis course might present the real numbers as the only complete ordered field.

The following nine axioms define a mathematical structure called a field.  Familar examples are $\R$ and the rational numbers ($\Q$ for short) .

\begin{enumerate}
\item {\em (commutativity of addition)} $a+b=b+a$
\item {\em (associativity of addition)} $a+(b+c)=(a+b)+c$
\item {\em (additive identity)} There is a number 0 such that $a+0=a$
\item {\em (additive inverses)} There is a unique number $-a$ such that $a+(-a)=0$
\\
\item {\em (commutativity of multiplication)} $ab=ba$
\item {\em (associativity of multiplication)} $a(bc)=(ab)c$
\item {\em (multiplicative identity)} There is a number 1 such that $1a=a$
\item {\em (multiplicative inverses)} There is a unique number $1/a$ such that $a(1/a)=1$
\\
\item {\em (distributive law)} $a(b+c)=ab+ac$
\end{enumerate}

Note that the integers are not a field, because there are no multiplicative inverses: $1/2$ is not an integer.

The two adjectives that distinguish the real numbers from the other fields are {\em complete}, meaning that there aren't any missing real numbers, unlike how $\sqrt2$ is missing from the rationals, and {\em ordered}, menaing that the symbols $<$ and $>$ are meaningful for $\R$; if $a$ and $b$ are distinct real numbers, then either $a<b$ or $a>b$.

We are more interested in the algebraic properties of the complex numbers, so we won't bother with the technical details of completeness and order.

\begin{exercise}Draw pictures of field axioms (1)-(9) for $\R$.\end{exercise}

\subsection{Complex arithmetic}

We will describe the complex numbers ($\C$ for short) geometrically, by defining a way to add or multiply points on a plane.  The procedure is similar to, and is in fact an extension of, the procedure for real numbers.  To add two points on the plane, we need to fix a point, which we'll call zero.  Once we have zero, we add complex numbers by stacking them next to each other.  But now the direction is important!  The procedure is the same as vector additon.  If $a$ and $b$ are points on the plane, thought of as complex numbers, and you draw arrows from 0 to $a$ and 0 to $b$, then $a+b$ is the number that you get by moving the arrow from 0 to $b$ so that the tail is at the head of the arrow for $a$.  If you do the same for $b+a$, you should get a parallelogram with vertices at 0, $a$, $b$, and $a+b$.

\begin{exercise}[Properties of Addition] The letters $a$, $b$, and $c$ stand for complex numbers.  Show by picture that this rule for addition satisfies the properties of addition in the field axioms.

\begin{enumerate}

\item Show that there is a number $-a$ such that $a+(-a)=0$. (Identity)

\item Show that $a+b=b+a$.  (Commutativity)

\item Show that $(a+b)+c=a+(b+c)$. (Associativity)

\end{enumerate}

\end{exercise}

Draw a horizontal line through 0, and choose a point to the right of 0 to call 1.  Notice that by repeatedly adding or subtracting 1 this gives us all of the integers as equally spaced points along a horizontal line in the complex number plane.  If we draw in the line connecting these dots, it corresponds to the real number line.

We know what it looks like to multiply real numbers, so we want to extend this idea that already works for the line to the whole plane.  This means that to find $ab$, while keeping 0 fixed, we stretch the  geometric thing that underlies the number system so that 1 is at $a$ and look at where $b$ goes.  Another way to say this is that you take the triangle with corners at 0, 1, and $b$, and you draw a similar triangle (one with the same angles in the same order) with the side similar to the one between 0 and 1 running between 0 and $a$.

\begin{exercise}[Properties of Multiplication]

The letters $a$, $b$, and $c$ stand for complex numbers.  

\begin{enumerate}

\item Show that there is a number $1/a$ such that $a(1/a)=1$. (Identity)

\item Show that $ab=ba$.  (Commutativity)

\item Show that $(ab)c=a(bc)$. (Associativity)

\item Show that the complex numbers have the distributive property, $a(b+c)=ab+ac$.

\end{enumerate}

\end{exercise}

Complex numbers do a lot algebraically that the real numbers can't do.  The most important thing about complex numbers is that negative numbers can have square roots! 

\begin{exercise} 
\begin{enumerate}

\item Draw a circle around 0 that passes through 1.  The number that is one quarter of the way around the circle, dierctly above 0 (on the perpendicular to the real line) is called $i$.  What happens when you multiply $i$ by itself?  (What is $i^2$?)

\item If you multiply any two numbers on this unit circle, what can you say about the result?

\item In any number system, $1^3=1$.  In the complex numbers, can you find a number other than 1 that when you take the third power you get 1?  Can you find another one?  They are on the unit circle.

\item How about fourth roots of 1?  (Hint: You already found one of them.)  

\item By now, maybe you can guess how to find $n$ different solutions to the equation $z^n=1$ for any positive integer $n$. 

\end{enumerate}
\end{exercise}

\subsection{Cartesian and polar forms}

This number $i$ is very special.  Lots of times, complex numbers are written in the form
$$
z=x+iy.
$$
Here, $x$ and $y$ are real numbers.  Starting from 0, $x$ tells us how far to go out horizontally, and $y$ tells us how far up to go vertically to find $z$ on the plane.  That is, $z$ is the point with Cartesian coordinates $(x,y)$ if our coordiantes on the plane put the origin at 0, $(1,0)$ at 1, and $(0,1)$ at $i$.  If $y=0$, then $z$ is a real number.  All of the surprising algebraic properties of $\C$ come from this $i$, this square root of $-1$.  Historically, taking square roots of negative numbers was rather hard to swallow, so $i$ or any multiple are called {\em imaginary numbers}.  Thus, we call $x$ and $y$ the {\em real} and {\em imaginary parts} of $z$, respectively, and often write $\Re z =x$ and $\Im z =y$.

\begin{exercise}
\begin{enumerate}

\item Draw the points $1+2i$, $1-2i$, $\frac{1+i}{\sqrt 2}$.

\item If $z=x+iy$, what is $-z$?

\end{enumerate}
\end{exercise}

We can use this representation and the distributive law to multiply complex numbers.  
$$
(x+iy)(s+it)=xs+iys+ixt+i^2yt=(xs-yt)+i(ys+xt)
$$
To divide complex numbers, observe that 
$$\begin{array}{rcl}
z\bar z &=& (x+iy)(x-iy)\\
            &=& x^2+y^2\\
            &=& |z|^2
\end{array}
$$
So $z\bar z/|z|^2=1$, or in other words $\bar z/|z|^2=1/z$, and if we want to compute $z/w=(x+iy)/(s+it)$, we can do this by computing 
$$
\frac{z\bar w}{|w|^2}=\frac{(x+iy)(s-it)}{s^2+t^2}=\frac{(xs+yt)+i(ys-xt)}{s^2+t^2}.
$$

\begin{exercise}
\begin{enumerate}

\item Calculate $\left(\frac{1+i}{\sqrt 2}\right)^2$, $(1+2i)(1-2i)$.

\item Calculate $\frac{2-i}{2+3i}$.

\item Find the $x$ and $y$ for the cube roots of 1 that you found above.

\end{enumerate}
\end{exercise}

The field of complex numbers is complete, but that isn't terribly important for us at this point.  It is important to realize, however, that the complex numbers are not ordered.  Which is greater, 1 or $i$?  The question has no answer.  The best we can do to compare two complex numbers is to give the {\em absolute value}, also called the {\em norm} or {\em modulus} or {\em magnitude}.  This is the distance from 0.  Just like for real numbers, the absolute value is denoted with vertical bars, and as it should, the complex notion of absolute value coincides with the real notion for real numbers inside the complex plane.   The absolute value can be calculated using the Pythagorean theorem if our number is written as $z=x+iy$: $|z|=\sqrt{x^2+y^2}$. 

There is another way to specify a point of the plane using coordinates.  Polar coordinates specify a point by giving the angle off the horizontal axis, sometimes called the {\em argument} and written $\arg(z)$, and the magnitude.  For a complex number, this representation is very useful, especially when multiplying complex numbers.

\begin{exercise}
If $z=(r,\theta)$ and $w=(s,\psi)$, what is $zw$ in polar coordinates?
\end{exercise}

The conversion between Cartesian ($z=x+iy$) and polar $z=(r,\theta)$ is straightforward.  To go from polar to Cartesian
$$
x=r\cos\theta, \quad y=r\sin\theta,
$$
and to go from Cartesian to polar
$$
r=|z|=\sqrt{x^2+y^2}, \quad \tan\theta=y/x.
$$

So we can write any complex number as $z=r\cos\theta+ir\sin\theta$.  If we factor out the $r$, we have $\cos\theta+i\sin\theta$.  Maybe you have seen exponentials and power series representations of functions.  If you haven't, just think of the following as a convenient notation, and a reason to be interested in learning about these things when they come up.  The power series expansions of sine and cosine are:
$$\begin{array}{rcl}
\sin(t)&=&t-\frac{t^3}{3!}+\frac{t^5}{5!}-\frac{t^7}{7!}+\cdots\\
\cos(t)&=&1-\frac{t^2}{2!}+\frac{t^4}{4!}-\frac{t^6}{6!}+\cdots.
\end{array}
$$
and the power series expansion of the exponential is 
$$
e^t=1+t+\frac{t^2}{2!}+\frac{t^3}{3!}+\frac{t^4}{4!}+\frac{t^5}{5!}+\cdots.
$$
Above, you found that $i$ is a 4th root of 1.  In particular: $i^0=1$, $i^1=i$, $i^2=-1$, $i^3=-i$, $i^4=1$, and then $i^5=i$ and the pattern repeats.  So if we write down the power series for $e^{i\theta}$ we get something interesting:
$$\begin{array}{rcl}
e^{i\theta}&=&1+i\theta+\frac{(i\theta)^2}{2!}+\frac{(i\theta)^3}{3!}+\frac{(i\theta)^4}{4!}+\frac{(i\theta)^5}{5!}+\cdots\\
&=&1+i\theta-\frac{\theta^2}{2!}-i\frac{\theta^3}{3!}+\frac{\theta^4}{4!}+i\frac{\theta^5}{5!}+\cdots\\
&=&\cos\theta+i\sin\theta.
\end{array}
$$
So we can write a complex number in polar coordinates as:
$$
z=re^{i\theta}.
$$
Famously, this expression gives rise to the equation $e^{i\pi}+1=0$.
\begin{exercise}
\begin{enumerate}
\item Convert $a=\frac{1}{\sqrt2}(1+i)$ and $b=1+i\sqrt3$ to polar coordinates and compute the product $ab$.
\item What are the $n^{th}$ roots of 1 in polar coordinates? 
\end{enumerate}
\end{exercise} 

When we defined $i$, we made an arbitrary choice.  If we had instead chosen the number one quarter of the way around the unit circle from one in the {\em clockwise} direction, we would have also found a number that squares to $-1$.  In algebraic terms, this is reflected in the fact that $(-i)^2=-1$.  The arbitrariness of this choice is reflected in a very important symmetry of the complex plane, called {\em complex conjugation}.  If $z=x+iy$ then write $\bar z=x-iy$, and call $\bar z$ the complex conjugate of $z$.

\subsection{Complex algebra}

Strictly speaking, this next part of the course isn't needed to understand the Hopf fibration, but I would feel deficient if I introduced the complex numbers and didn't talk about these following ideas.

Complex numbers give us the abilily to solve algebraic equations.  The $n^{th}$ roots of 1 are the solutions to the equation $x^n=1$.  If our variable $x$ can take complex values, then we can find $n$ roots for {\em any} polynomial of degree $n$.  This result is so important that it gets a name signifying how powerful it is.

\begin{thm}[The Fundamental Theorem of Algebra]
Any polynomial of degree $n$ with coefficients in $\R$ (or even $\C$) can be factored into linear terms.
\end{thm}

Proving this theorem rigorously would take us too far afield for now.  There are many different ways that it can be proved.  Perhaps later in your mathematical development, you will get to decide which ones are your favorites.  Some parts of my favorite proof will be described below.

\begin{exercise}
\begin{enumerate}
\item Find the roots of $z^2+4z+3=0$.
\item Show that if $p(z)$ is a quadratic polynomial with real coefficients and $z_1$ is a root, then $\bar z_1$ is also a root.
\item Graph the parabola defined by $y=x^2+4x+5$ in the plane $\R^2$.  Revise the description of the three types of parabolas that relies on the discriminant in the quadratic formula to one based on complex numbers.
\item Find the roots of $z^3-3z^2+4z-2$.
\end{enumerate}
\end{exercise}

\subsection{Functions of a complex variable}

Perhaps the most important things to study about a number system are the functions.  Most of the functions of a real variable that you are familiar with also make sense when the input and output are thought of as complex.  For example, the function $f(x)=x^2$ makes perfect sense for $x$ a complex number.  

There are actually a lot of very significant differences in the theory of functions of a real variable and the theory of functions of a complex variable, but the one that we'll pay the most attention to is the simple fact that because $\C$ is a two dimensional space when viewed with our ``real eyes'' the notion of drawing a graph of a function, as we do with $f(x)=x^2$ when we draw a pair of axes and a parabola passing through the point where they cross, is simply impossible in a three dimensional space.  We would need two dimensions for the input, and two dimensions for the output, for a total of four.  

All is not lost, however, and we can graphically visualize complex functions by the ways that they transform shapes drawn on the plane.  For example, multiplying by a complex number $z$ rotates by $\arg(z)$ and scales by $|z|$, and complex conjugation reflects in the real axis.  For more complicated functions, we can develop an visual understanding by looking at how a grid is transformed.   

\begin{exercise} This exercise shows how the square grid is transformed by the function $f(z)=z^2$.  
\begin{enumerate}
\item Find the real and imaginary parts of $z^2$ if $z=x+iy$.
\item If $w=u+iv=z^2$, describe the shape in the $w$-plane that is the image of the square with sides defined by the lines $x=0$, $x=1$, $y=0$ and $y=1$.
\item Do the same for the similar squares in the second third and fourth quadrants.
\item Do the same for the similar squares of twice the side length and half the side length in the first in the first quadrant.
\item Describe the action of the function $z\mapsto z^2$.   Include in this description some reference to why the graph in $\R^2$ of $x\mapsto x^2$ looks the way it does.
\item Now look at $z\mapsto z^3$, and describe the transformation caused by this function.
\item How about $z\mapsto z(z-2)$?
\end{enumerate}
\end{exercise}

We now introduce an important tool in the study of complex functions.

\begin{defn} A {\em path} in $\C$ is a continuous map $C:[0,1]\to \C$.  A path is called {\em closed} if $C(0)=C(1)$.\end{defn}

\begin{defn} The {\em winding number} of a closed path $C$ is defined as the number of times the path moves around 0 counter-clockwise.\end{defn}

If $f:\C\to\C$ is a continuous function, we can learn a lot about it by looking at the winding numbers of various paths $f^{-1}(C)$, where $C$ is a closed path.

\begin{exercise} Let $C$ be the path tracing out the unit circle: $C(t)=e^{2\pi i t}$, and let $K$ be the path tracing out a circle of radius one around the number 2: $K(t)=e^{2\pi i t}+2$.
\begin{enumerate}
\item What are the winding numbers of $C$ and $K$?
\item If $f(z)=z^2$, what are the winding numbers of $f^{-1}(C)$ and $f^{-1}(K)$?
\item What if $f(z)=z^n$?
\item $f(z)=z(z-2)$?
\end{enumerate}
\end{exercise}

What these examples are getting at is the fact that the winding number can be used to detect zeros of a complex function.  This can be used to prove the Fundamental Theorem of Algebra in the following way.  For $z$ with $|z|>>0$, any polynomial of degree $n$ looks enough like  $z^n$ that a circle with this large radius will have winding number $n$, so there are $n$ zeros inside the circle.  Making this precise requires some careful work, so that's all we'll say in that direction.

\newpage

%% file: spheres.tex
\section{Spheres}

You might think of the sphere as the set of points defined by the equation
$$
x_1^2+x_2^2+x_3^2=1.
$$
This defines a surface in space, consisting of those points at distance one from the origin.

Similarly, you might think of a circle as the solutions to the equation 
$$
x_1^2+x_2^2=1,
$$
as this defines a curve on the plane, consisting of those points at distance one from the origin.

What about the solutions to the even simpler equation
$$
x_1^2=1,
$$
or the more complicated
$$
x_1^2+x_2^2+x_3^2+x_4^2=1?
$$
It makes sense to mathematicians to call all of these objects spheres. The dimension of a sphere is the number of parameters required to specify a point.  So we say that the circle is a one dimensional sphere, or one sphere for short, or $S^1$ for even shorter, that the surface of the earth is a two sphere, or $S^2$ (approximately - the Earth isn't exactly round, but it is pretty close), and by analogy, the set of points in four dimensional space satisfying the equation $x_1^2+x_2^2+x_3^2+x_4^2=1$ is a three dimensional sphere or $S^3$.   

\subsection{Dimension} ``Wait a minute!'' you might say, ``Four dimensional space, how the heck am I supposed to imagine that?!''  Or, maybe you have thought about it a bit, and have a few ideas.  Anticipating this, one mathematician I know will occasionally begin a talk about dimensions with the rhetorical question, ``So, is the fourth dimension time, or what?''

To a mathematician, this question need not be any more meaningful than the reply, ``No, the second dimension is time, the fourth dimension is red.''  Mathematically, a dimension is something that can be measured, something that can take a value.  It is a characteristic of an object that in some rough way, describes its complexity.  A phrase in common usage that aligns closely to the mathematician's notion of dimension is, ``that adds a new dimension to the situation.''

A circle is a relatively simple object in that a single number, for example the angle measured counterclockwise from a fixed point, is enough to fully specify any point on the circle.  One measurement locates a point on the circle, and we call it a one dimensional object.  On the other hand, the weather is a very high dimensional system.  It is true that the temperature in Seattle today, the time of year, what the winds and clouds are like over the Olympic peninsula, and what the weather will be like tomorrow are all related, but the relationship is very complex and depends on a large nuimber of variables.  So when a meteorologist makes a prediction of the weather based on all the information at hand, she is using some sort of model of the weather that takes all of these measurements into account.  We say that this model is a high dimensional model.

For objects with a small number of dimensions, it can be helpful to use our familarity with $\R^3$, coming from the fact that our surroundings look very much like this object, to visualize these spaces.  It turns out that there are many three dimensional objects that we are able to ``see'' with a little bit of imagination and effort.

Four dimensions are not too hard to visualize.  We can think of time as a fourth dimension; a solid four dimensional thing is something that exists for some span of time as a solid three dimensional thing.  But there is no need to insist on using time.  Another popular choice is color.  

\begin{exercise} Using color as a fourth dimension, explain why you can't tie a knot in a piece of string in a four dimensional space.
\end{exercise}

However, directly visualizing things in four dimensions is generally only good enough to let us see topological facts.  If we try to visualize the difference between a nice round three dimensional sphere in $\R^4$ and some topologically identical distortions that are not so round, it can be less than straightforward.  Therefore, instead of trying to build a globe of the three sphere, we can make a map.  This map will be drawn on a three dimensional Euclidean space, just as a map of the two dimensional sphere such as the earth is drawn on a two dimensional Euclidean space such as a piece of paper.  But before we even talk about the three dimensional sphere at all, we will look at some properties of the lower dimensional spheres, so that we will know what to look for.

\subsection{$S^0$: the zero sphere}

A zero dimensional sphere is a pair of points.  There are two solutions to the equation $x_1^2=1$, namely $x_1=1$ or $x_1=-1$.  So we can write $S^0=\{-1,1\}$.  There isn't much more to it than that.  It is the only sphere that is disconnected, in that the two points are separated by some distance, but other than that, it isn't terribly interesting.  On the other hand, it shows up in its proper place whenever we need it to.

\subsection{$S^1$: the circle}

A one dimensional sphere is a circle, the set of points in the plane equidistant from a fixed center, or the solutions to the equation $x_1^2+x_2^2=1$.  There are a couple of special ways in which it comes up that are worth mentioning.

The complex numbers $z$ with $|z|=1$ are called the {\em unit complex numbers} or if it is understoood that we are talking about complex numbers, just the {\em units}.  The units form a circle of radius one.  One crucial feature of the units is that the product of two units is again a unit.  Just as numbers as measurement of distances or lengths are represented graphically as the line, we can use the circle to represent numbers as measurements of rotations.  In polar coordinates on $\C$, any point on the circle can be written $z=e^{i\theta}$.  It is often handy to have such a compact notation.

\subsubsection{Topology} When you ``circle back around to have another look,'' or lament that you are ``moving in circles,'' it is not that you mean to say that you have traced an arc of points equidistant from a central point, but rather that your path has carried you back to your starting point.  To a topologist (which is a special type of mathematician, quite similar to a geometer), the roundness of a circle (or other spheres for that matter) is not intrinsic.  In the idiom of topology, it makes sense to take an extension cord, plug one end into the other, and call it $S^1$, even if it is not even close to being shaped into a circle.  Topology looks at properties of geometric things that are intrinsic in a way that the particular shape it happens to take is not.  Since the extension cord plugged into itself could be arranged into a perfect circle without cutting it up and rearranging it (we might have to untie some knots by unplugging it, untying the knot, and plugging it back together, but since it gets put back together exactly the way it was before, that's okay), it is in a topological sense, the same as a circle. 

\subsubsection{The projective line}  Thinking topologically about the circle allows us to see it as another interesting geometric object.  The {\em real projective line} is the set of lines through the origin in the plane.  Since each line is designated by it's slope, this is also the set of ratios $[x:y]$, where $x$ and $y$ real numbers.  (We say that a vertical line has infinite slope; it corresponds to the ratio $[0:1]$.)  We write $\Pp_\R^1$ for the projective line.  

If we think of the unit circle in the plane, each line through the origin makes a diameter of that circle.  So in particular, if we take the upper semicircle $\{(x,y)|x^2+y^2=1, y\geq0\}$, there is one point on the semicircle for each line, except for the two points $(1,0)$ and $(-1,0)$, which both correspond to the same line: the $x$-axis.  So if we identify these two points, as if the semicircle is an extension cord and we plug the ends together, we get a circle, topologically.

We can also represent (most of) the projective line as a straight line.  If we draw the line $x=1$ in the plane, each line intersects this line at the point where $y$ equals it's slope.  The vertical line is missing from this representation, but it is the only one, so we can say that this line $x=1$ gets an extra point, called the {\em point at infinity}, and that the projective line is thus $\R$ plus one more point, called $\infty$.

That the projective line is circle can also be seen via an operation called stereographic projection, which is an important way to visualize spheres that works in all dimensions.  Stereographic projection gives a way to map an $n$-dimensional sphere on an $n$-dimensional Euclidean space.  Consider the configuartion already studied, with the line $x=1$ giving a map of $\Pp_\R^1$, and now draw a circle with diameter the unit interval on the $x$-axis.  For any point $p$ on the circle, other than the origin, the line through the origin and $p$ hits both the circle and the line $x=1$ at exactly one point.  So the circle is the projective line; {\bf$\infty$}, or $[0:1]$, is the origin, and the remainder of the points $p=(p_1,p_2)$ on the circle are represented by the ratio $[p_1:p_2]$, or the point $(1,p_2/p_1)$ on the line $x=1$.

\subsubsection{Projective transformations}

One of the fundamental features of any geometric object that can help us to understand it are the transformations that can be used to move it around without fundamentally changing what it is.   The projective line $\Pp^1_\R$ is the collection of lines through the origin in $\R^2$, so any transformation of $\R^2$ that sends lines through the origin to lines through the origin is also a transformation of the projective line.  In other words, any linear transformation of the plane with nonzero determinant gives rise to a transformation of the projective line.  These transformations correspond to choosing various different lines in the plane to project onto.

\begin{exercise} None of the following characteristics are preserved by projective transformations of the line.  Give examples to show this.
\begin{enumerate}
\item Any particular point: there is no ``origin'' of the projective line.
\item The distance between two points.
\item A point being between two other points.
\end{enumerate}
\end{exercise}

There is a quantity that is preserved if we consider four points on $\Pp^1_\R$.  This is called the {\em cross ratio}.   If $a$, $b$, $c$, and $d$ are any four distinct points on $\Pp^1_\R$, choose coordinates (by choosing a line to project onto and a scale to use on that line) and define the cross ratio as the fraction
$$
(a,b;c,d)=\frac{(a-c)(b-d)}{(a-d)(b-c)}.
$$
This is unchanged if we apply a transformation $T$.  

\begin{exercise} 
A projective transformation $T$ can be represented computationally by the matrix of a corresponding linear transformation of $\R^2$.  Show that such a $T$ does not change the cross ratio.
\end{exercise}

\begin{exercise}
The permutation group $S_4$ acts on the cross ratio by reordering  $a$, $b$, $c$, and $d$.  However, for a fixed  $a$, $b$, $c$, and $d$, there are only six different values that the cross ratio can take.  Suppose that $(a,b;c,d)=\lambda$.  
\begin{enumerate}
\item Which permutations of $a$, $b$, $c$, and $d$ leave the value $\lambda$ invariant?
\item For permutations that do change the value of $\lambda$, what are the resulting values?
\end{enumerate}
\end{exercise}

\subsection{$S^2$: the skin of the globe}

A two dimensional sphere is what most people think of when they think of a sphere.  It is the surface of a ball or a globe.  Its geometry is quite rich, and by exploring some aspects of it with analogs in higher dimensions, we will be prepared to understand the structure of $S^3$ via analogies with $S^2$.

\begin{defn} A {\em great circle} is a circle on the sphere that is as long as possible.  Great circles are the intersections of planes through the center of the sphere with the sphere.  The shortest distance between two points on the surface of the sphere is along a great circle.
\end{defn}

\begin{exercise}
\begin{enumerate}
\item Consider the triangle on the unit sphere $x_1^2+x_2^2+x_3^2=1$ made by the intersection of the coordinate planes $x_i=0, \quad i=1\ldots3$ with the sphere in the orthant where all the coordinates are non negative ($x_i\geq0, \quad i=1\ldots3$.  What are the angles of this triangle?  What is its area?  (The surface area of a sphere with radius $r$ is $4\pi r^2$.)
\item Find the area of a spherical triangle with angles $\alpha$, $\beta$, and $\gamma$ on the unit sphere.
\end{enumerate}
\end{exercise}

\subsubsection{The Riemann sphere and stereographic projection}

Like the circle, $S^2$ also represents a projective line, but it is a projective line for the complex numbers.  In this guise, it is also called the Riemann Sphere.  Analogous to $\Pp_\R^1$, the complex projective line $\Pp_\C^1$ is defined by ratios of complex numbers $[z_1:z_2]$, and if $z_1\neq0$, this can be thought of as the ``slope'' $z_2/z_1$, though how exactly this is a slope is hard to visualize, since $\C^2$ is four dimensional as a real space.  There is again a point which we call $\infty$, corresponding to the ratio $[0:z_2]$.

Stereographic projection gives us a way to see $\Pp_\C^1$ as $S^2$.  Since any point $[z_1:z_2]\in\Pp^1_\C$ other than $\infty$ can also be represented by $[1:\frac{z_2}{z_1}]$, we can define $z=x+iy=\frac{z_1}{z_2}$, and put this plane into $\R^3$ with coordinates $(x,y,t)$.  Now put a sphere into this $\R^3$ so that the unit interval on the $t$-axis is a diameter (so the center of the sphere is $(0,0,1/2)$, and it has radius 1/2).  We can project from the point $(0,0,1)$: for every point $p$ on the sphere except $(0,0,1)$, draw the line through $p$ and $(0,0,1)$.  This line hits the sphere at these two points, and it hits the plane $t=0$ exactly once.  

\begin{exercise}

\begin{enumerate}
\item Show that the line between $(p_1,p_2,p_3)$ and $(0,0,1)$ also passes through $(\frac{p_1}{1-p_3},\frac{p_2}{1-p_3},0)$.  So a formula for stereographic projection from the sphere to the plane is given by 
$$
\sigma(p_1,p_2,p_3)=\left(\frac{p_1}{1-p_3},\frac{p_2}{1-p_3}\right),
$$ 
\item Show that the inverse map from the plane to the sphere $x^2+y^2+(t-1/2)^2=1/4$ is defined by 
$$
\sigma^{-1}(x,y)=\left(\frac{x}{x^2+y^2+1},\frac{y}{x^2+y^2+1}, \frac{x^2+y^2}{x^2+y^2+1}\right)
$$
\item Write these formulas for $\sigma$ and $\sigma^{-1}$ in terms of complex numbers.
\item There is no particular need to use the sphere $x^2+y^2+(t-1/2)^2=1/4$.  In fact, show that the formula for $\sigma$ found in part (1) also defines stereographic projection from the sphere $x^2+y^2+t^2=1$ to the plane $t=0$.  Show that if we want to invert stereographic projection from the plane back onto the unit sphere, that we use the formula
$$
\sigma^{-1}(x,y)=\left(\frac{2x}{x^2+y^2+1},\frac{2y}{x^2+y^2+1}, \frac{x^2+y^2-1}{x^2+y^2+1}\right).
$$

\end{enumerate}
\end{exercise}

Via stereographic projection we can match each point on the sphere with a complex number, except for $(0,0,1)$.  So we define $\sigma(0,0,1)=\infty$.

\begin{exercise} 

\begin{enumerate}
\item What points on the sphere correspond to the unit complex numbers?
\item Consider the transformation $z\mapsto iz$.  Describe the transformation of the sphere that this corresponds to.
\item Consider the transformation $z\mapsto 2z$.  Describe the transformation of the sphere that this corresponds to.
\item Consider the transformation $z\mapsto1/z$.  What transformation of the sphere does this correspond to?
\item Describe the sequence of points 
$$
\left\{ \ldots, \frac1{(1+i)^3}, \frac1{(1+i)^2}, \frac1{1+i}, 1, 1+i, (1+i)^2, (1+i)^3, \ldots \right\}
$$ 
on the sphere.  There is an Escher print illustrating this sort of transformation.
\item Find the coordinates in $\R^3$ for a regular octahedron sitting inside the sphere with vertices at $\sigma^{-1}(0)$, $\sigma^{-1}(1)$, and $\sigma^{-1}(\infty)$.  Where are the rest of the vertices when projected to the plane?  If you connect the vertices with great circles on the surface of the sphere, what do these great circles map to via stereographic projection?
\item  Instead of an octahedron, put a cube inside the sphere with its faces parallel to the coordinate planes.  What are the complex numbers corresponding to the vertices of the cube?  If the vertices of the cube are connected with arcs of great circles instead of straight lines, what do the projections of these great circle edges of the cube look like?  

What if the cube was sitting with vertices at $\sigma^{-1}(0)$ and $\sigma^{-1}(\infty)$, and two more of its vertices along the circle corresponding to the real axis.  Where are the rest of the vertices, and where are the great circle edges?
\end{enumerate}
\end{exercise}

\subsubsection{M\"obius transformations}

There is a lot of interesting geometry to the complex numbers and the Riemann sphere that have to do with {\em linear fractional transformations}, also called {\em M\"obius transformations}.  These are functions of the form
$$
f(z)=\frac{az+b}{cz+d}
$$
where $a,b,c$ and $d$ are complex numbers.   M\"obius transformations naturally act on the Riemann sphere; even when $cz+d=0$, we would like to have a value for $f(z)$, and it makes sense to set this value as $\infty$.  Conversely, if we compute $f(z)$ for values of $z$ with $|z|$ growing larger and larger, the value of $f(z)$ approaches $a/b$.  So we can set $f(\infty)=a/b$.

\begin{exercise}
Show that a M\"obius transformation $f(z)=\frac{az+b}{cz+d}$ is determined by the values of any three points in two steps:
\begin{enumerate}
\item Show that for any three distinct inputs, $z_1, z_2, z_3\in\C$ there is a unique M\"obius transformation $h(z)$ such that $h(z_1)=0$, $h(z_2)=1$ and $h(z_3)=\infty$.
\item Show that for any three distinct values, $w_1, w_2, w_3\in\C$ there is a unique M\"obius transformation $g(z)$ such that $g(0)=w_1$, $g(1)=w_2$ and $g(\infty)=w_3$.
\end{enumerate}
So by setting $f(z)=g(h(z))$, we get a unique M\"obius transformation with $f(z_1)=w_1$, $f(z_2)=w_2$ and $f(z_3)=w_3$.
\end{exercise}

\begin{exercise}
Since the definition of the cross ratio as an invariant of the projective line didn't really rely on the fact that we were thinking about $\Pp^1_\R$, it is also an invariant of $\Pp^1_\C$.  Show that if $z_1$, $z_2$, $z_3$, and $z_4$ are distinct complex numbers, and $f(z)$ is a M\"obius transformation with $f(z_1)=0$, $f(z_2)=1$ and $f(z_3)=\infty$, then $f(z_4)=(z_1,z_2;z_3,z_4)$
\end{exercise}

One of the important things about M\"obius transformations is that they can be composed by matrix multiplication.  A M\"obius transformation
$$
f(z)=\frac{az+b}{cz+d}
$$
can be encoded in the matrix
$$
M_f=\twobytwo abcd.
$$
\begin{exercise}
Show that if 
$$
f_1(z)=\frac{a_1z+b_1}{c_1z+d_1}, \quad \mbox{and} \quad f_2(z)=\frac{a_2z+b_2}{c_2z+d_2}, 
$$
then the composition $f_1(f_2(z))$ is given by the matrix $M_{f_1}M_{f_2}$.
\end{exercise}

\subsubsection{Sphere inversion}

\begin{figure}
\centerline{\resizebox{0.8\textwidth}{!}{\includegraphics{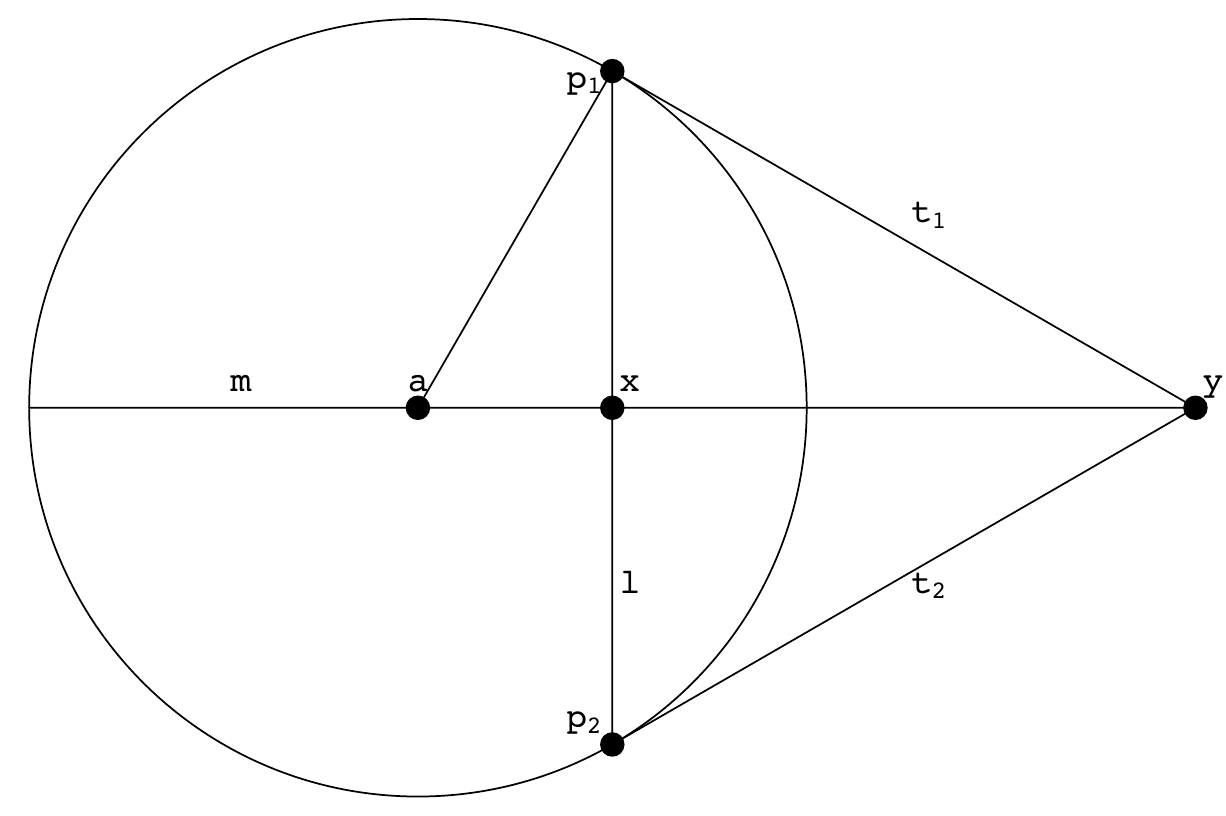}}}
\caption{A construction of inverse points}
\label{invert1}
\end{figure}

Sphere inversion is a way of turning space inside out, so that the inside of a sphere is sent to the outside, and vice versa.  It will help us to understand the geometry of stereographic projection.

Sphere inversion is defined for spheres of any dimension.  
\begin{defn} Two points $p$ and $q$ are said to be {\em inverse with respect to a sphere $S$ with center $a$ and radius $r$} if the following conditions are satisfied:
\begin{itemize}
\item $p$, $q$ and $a$ are colinear, and $a$ is not between $p$ and $q$.
\item $|ap||aq|=r^2$
\end{itemize}
\end{defn}
This is enough to completely specify $q$ given $p$, so for any sphere $S$ there is a transformation $i_S$ that takes each point $p$ to its inversion in $S$.  Note that the definition is symmetric, in that $p$ and $q$ can be interchanged, so $i_s(p)=q$ means than also $i_S(q)=p$.  The only caveat is that $i_s(a)$ is not defined.  To remedy this, we can add a point at infinity, and say that $i_S(a)=\infty$.  This is reminiscent of stereographic projection, and for good reason, as we will see.

\begin{exercise} Show that in the configuration of figure \ref{invert1} the points $x$ and $y$ are inverse.  The figure represents the following construction: to invert a point $x$ which lies inside a circle $C$ with center $a$, draw a line $l$ through $x$ which is perpendicular to the line $m$ containing $x$ and $a$.  Line $l$ will meet $C$ at two points, call them $p_{1}$ and $p_{2}$.  Draw lines $t_{1}$ and $t_{2}$,
tangent to $C$ at $p_{1}$ and $p_{2}$, respectively.  The intersection
point $y$ of $t_{1}$ and $t_{2}$ lies on $m$, by symmetry, and we say
that $y$ is inverse to $x$ with respect to $C$.  This construction can
easily be reversed for a point outside of $C$, given $y$, draw $t_{1}$
and $t_{2}$, connect them with line $l$, and the inverse point $x$ is
the intersection of $l$ and $m$, $m$ being drawn in this case by
connecting $y$ to $a$.

\end{exercise}

One of the most important properties of inversion is that the inversion of a sphere is a sphere.  That is, if $K$ is a sphere of any dimension up to and including the dimension of $S$, then $i_S(K)$ is also a sphere.  The only caveat here is that if $K$ passes through the center of $S$, then $i_S(K)$ is a plane, which we can think of as a sphere containing $\infty$.

\begin{figure}
  \centerline{\resizebox{0.5\textwidth}{!}{\includegraphics{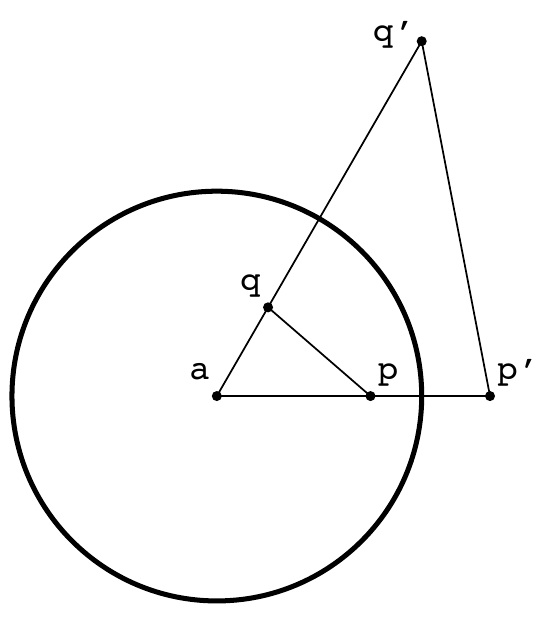}}}

  \caption{In this configuration $p$ is inverse to $p'$ and $q$ is inverse to $q'$.}
\label{invert7}
\end{figure} 

\begin{figure}
  \centerline{\resizebox{0.5\textwidth}{!}{\includegraphics{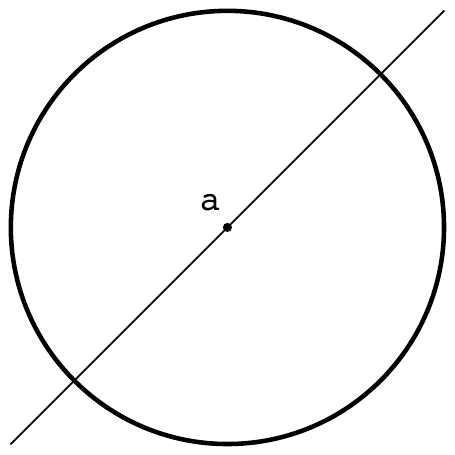}}}
	\caption{$a \in \ell$}
\label{invert2}
\end{figure}

\begin{figure}
  \centerline{\resizebox{0.8\textwidth}{!}{\includegraphics{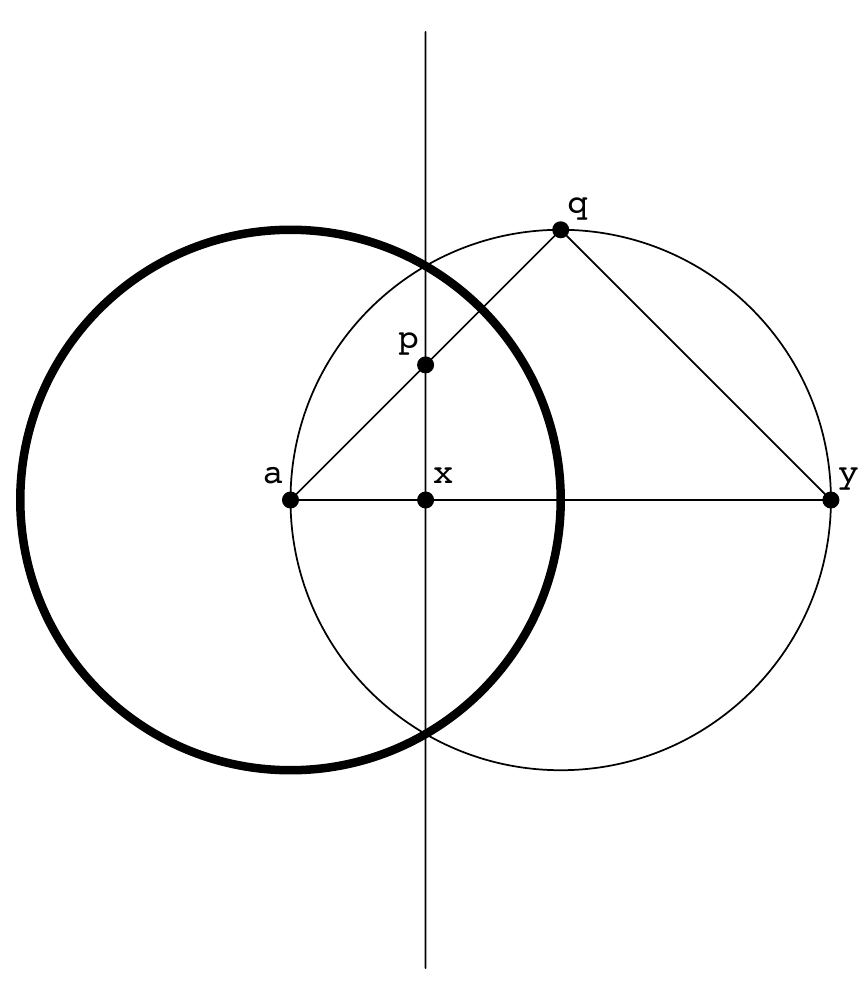}}}
	\caption{$a \not \in \ell$}
\label{invert3}
\end{figure}

\begin{figure}
   \centerline{\resizebox{0.8\textwidth}{!}{\includegraphics{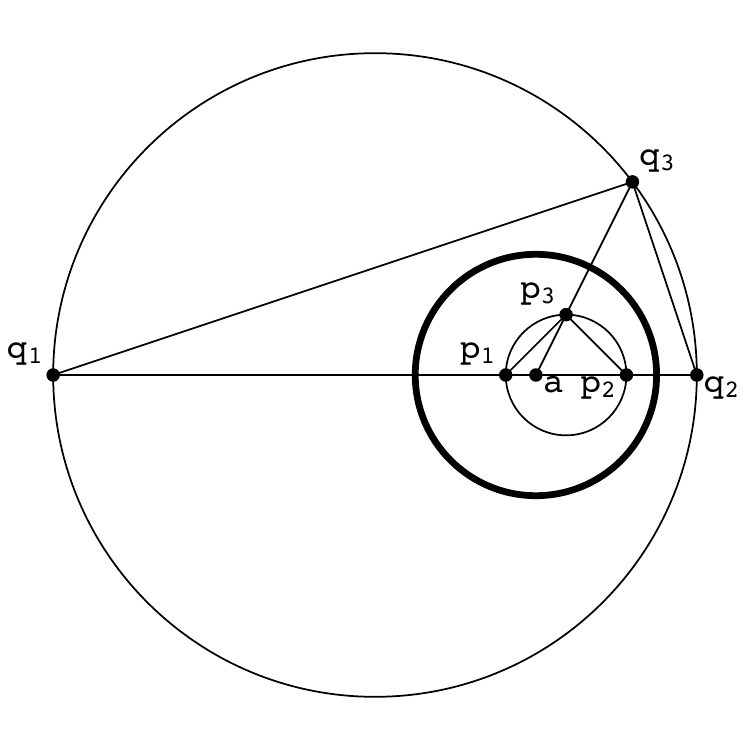}}}
	\caption{$a$ inside $K$}
\label{invert4}
\end{figure}

\begin{figure}
  \centerline{\resizebox{0.8\textwidth}{!}{\includegraphics{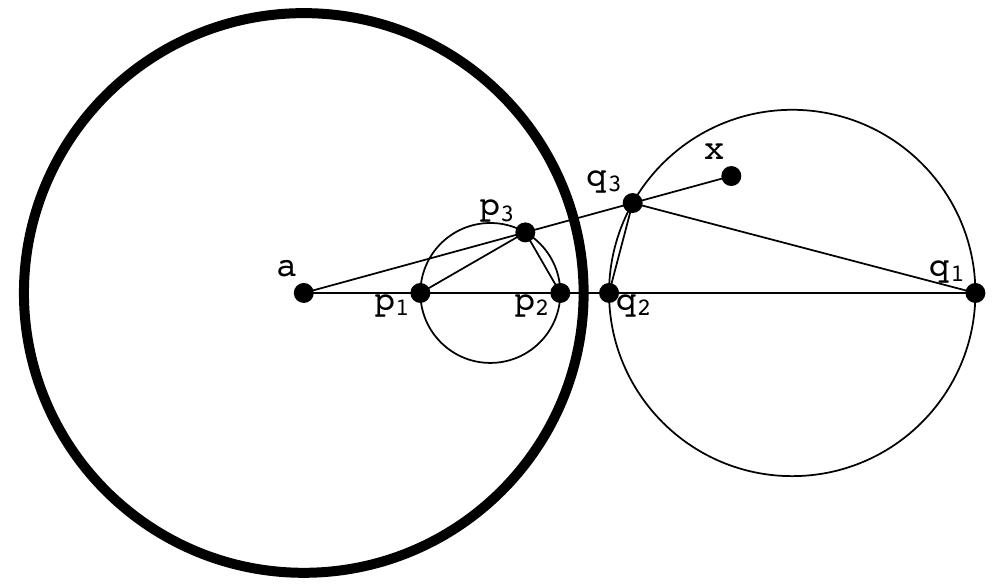}}}
	\caption{$a$ outside $K$}
\label{invert5}
\end{figure}

\begin{exercise} This exercise shows that the inversion in $\R^2$ of a circle is or a line is a circle or a line.  (Circles passing through the center of the inverting circle are sent to lines, and vice versa.)  There are a few cases to check.  Most of them rely on the so called {\em Carpenter's angle theorem}: a triangle inscribed a circle such that one side is a diameter, is a right triangle with the diameter as hypotenuse.

In what follows, $C$ is a circle with center $a$, shown as the darker circle in the figures.
\begin{enumerate}
\item Show that in the configuration of figure \ref{invert7}, $ \angle aqp = \angle ap'q' $ and $ \angle apq = \angle aq'p' $.

\item If $\ell$ is a line passing through $a$ as in figure \ref{invert2}, show that $i_C(\ell)=\ell$.

\item If $\ell$ is a line which does not pass through $a$, use  figure \ref{invert3} to show that the inverse $I_C(\ell)$ is a circle containing $a$.

\item If $a$ is inside a circle $K$, use figure \ref{invert4} to show that $i_C(K)$ is a circle.

\item If $a$ is outside a circle $K$, use figure \ref{invert5} to show that $i_C(K)$ is a circle.

\end{enumerate}
\end{exercise}

We can also explicitly write down the formula for inversion of a point $z$ in a circle $C$ with center $a$.
$$
i_{C}(z)=a + \left(\frac{r}{|z-a|}\right)^{2}(z-a)
$$ 
Note that this formula is not particular to any dimension, it works for all spheres, not just circles.
\begin{exercise}
\begin{enumerate}
\item Rewrite the formula for $i_C(z)$ in terms of complex conjugation.
\item Compose two inversions, that is, take $i_C(i_K(z))$.  Show that this is a M\"obius transformation.
\end{enumerate}
\end{exercise}

One of the most important features of stereographic projection is that circles on the sphere correspond to circles on the plane, with the exception that circles on the sphere that pass through the point of projection (generally $(0,0,1)$ in these notes) correspond to straight lines in the plane.  

\begin{exercise}  Find a sphere $R$ such that $i_R$ agrees with $\sigma$ when restricted to the sphere.
\end{exercise}

The existence of this way to do stereographic projection allows us to prove one of the most important facts about stereographic projection.

\begin{thm}\label{stereocircles}
If $C$ is a circle on a sphere $S$, and $\sigma$ is stereographic projection from a point $p\in S$ to a plane $M$, then $\sigma(C)$ is a circle on $M$, unless $p\in C$, in which case $\sigma(C)$ is a straight line.  Conversely, if $K$ is a circle or a line on $M$, then $\sigma^{-1}(K)$ is a circle on $S$.
\end{thm}

\subsubsection{Circles of Apollonius}

\begin{figure}
\begin{tabular}{c}\resizebox{\textwidth}{!}{\includegraphics{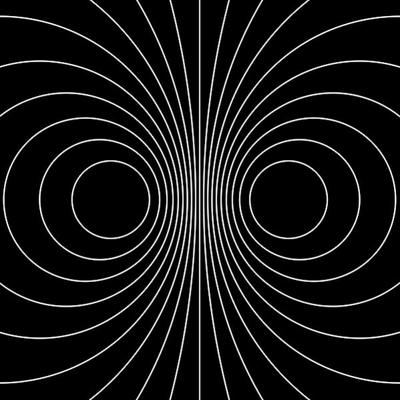}\qquad \qquad \includegraphics{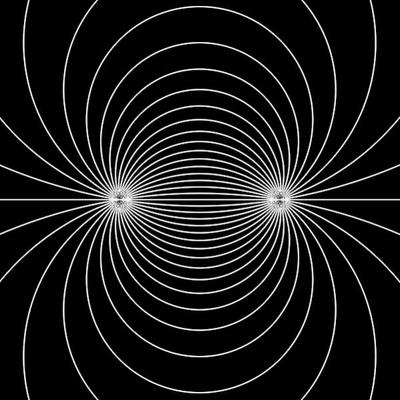}} \\ \\  \resizebox{0.6\textwidth}{!}{\includegraphics{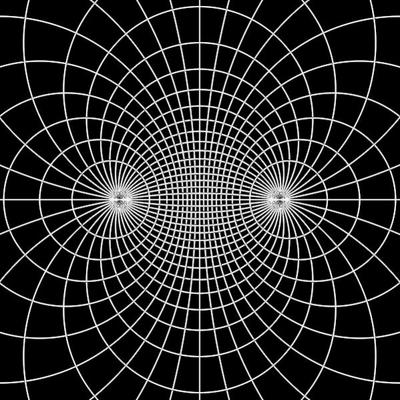}}\end{tabular}
\caption{Apollonian circles: the hyperbolic family {\em (left)}, the elliptic family {\em (right)} and the full configuration {\em (bottom)}}
\label{apollonius}
\end{figure}

There is a configuration of circles that has fascinated geometers for millennia, called the {\em circles of Apollonius}, that is intimately tied to stereographic projection and inversions.  The circles of Apollonius, shown in figure \ref{apollonius}, depend on a pair of distinct points in a plane, we'll call them $p$ and $p'$, and consist of two families of circles.  Each family having one degenerate member that is a straight line.  One family consists of all circles passing through both $p$ and $p'$.  This is called the {\em elliptic} family.  The second family consists of all circles such that the inverse of $p$ is $p'$, and this family is called the {\em hyperbolic} family

When two circles intersect, we define the angle of intersection as the angle between the tangent lines at the point of intersection that is on the outside of both circles.  If this is a right angle, then we say that the circles are {\em orthogonal}.  

\begin{exercise}
\begin{enumerate}
\item Show that a circle $K$ is orthogonal to a circle $C$ if and only if $i_C(K)=K$.
\item Show that for circles of Apollonius, any circle in the elliptic family is orthogonal to any circle in the hyperbolic family.
\end{enumerate}
\end{exercise}

This exercise implies that the full configuration is symmetric with respect to inversion in any one of the constituent circles.

To connect the circles of Apollonius to stereographic projection, start with two antipodal points on the sphere, such as the north and south poles.  Antipodal points determine a family of great circles, those that pass through both points, like the meridians of longitude on the earth.  The two points also determine another family of circles, those that have those two points as their centers, like the parallels of latitude.  The stereographic projection of these two families of circles on the sphere (from a point other than one of the two antipodal points, so perhaps you'll need to stereographically project from magnetic north) are the two families that make up the circles of Apollonius on the plane.  The meridians create an elliptic family, and the parallels form a hyperbolic family.

\begin{exercise} \label{describegreatcircles} Describe the images of the great circles on the sphere under stereographic projection in terms of how they intersect the unit circle on the plane.
\end{exercise}

\subsection{$T^2$:  the torus}

\begin{figure}
  \centerline{\resizebox{\textwidth}{!}{\includegraphics{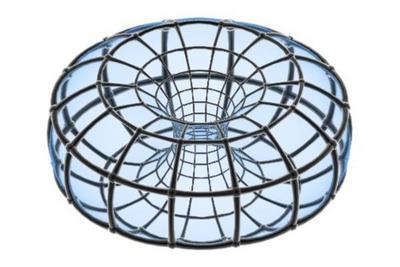}}}
	\caption{A torus}
\label{torus}
\end{figure}

The torus ({\em pl. tori}) is not a sphere, as even a topologist could tell you.  However, it is an important surface to be familiar with in the study of the three dimensional sphere so we'll take a look at it.  The torus (colloquially the doughnut, and written $T^2$ for short) is the mathematical name for the surface of an inner tube.  Its essential characteristic is that there are two circles making up the surface; the one that makes it roll, and the one that keeps the air in.

There are tori of every dimension, just as there are spheres and Euclidean spaces of each dimension.  Just as a point in $\R^n$ can be specified by an $n$-tuple of points from $\R$, a point on the $n$-torus $T^n$ is specified by an $n$-tuple of points from the circle $S^1$.  So a point on a two dimensional torus is specified by a pair $(\theta,\psi)\in S^1\times S^1$.

A {\em torus of revolution} is a torus in $\R^3$ that is defined by rotating a circle $C$ about an axis $\ell$ that lies in the plane of $C$ but does not intersect $C$.  An inner tube is a torus of revolution - the axis is the axle of the wheel, and the circle $C$ is a radial slice of the tube.  Any torus of revolution can be described by two numbers: the large radius, or the distance from the axis of revolution to the center of $C$, and the small radius - the radius of $C$.

\begin{figure}
  \centerline{\resizebox{0.8\textwidth}{!}{\includegraphics{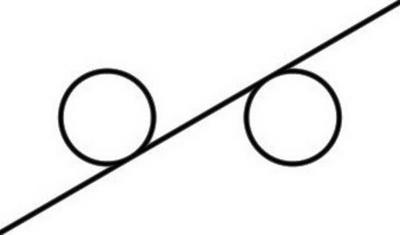}}}
	\caption{A cross section of a plane cutting a torus of revolution}
\label{hopfcut}
\end{figure}

\begin{exercise}  \label{toruseqns}
\begin{enumerate}
\item Starting with the parametric form for a circle $(\cos t, \sin t)$, find a parametric form for a torus of revolution with large radius $R$ and small radius $r$.  Use $\theta$ and $\psi$ for the parameters of your torus, so you are looking for a triple of functions  $(x_1(\theta,\psi), x_2(\theta,\psi), x_3(\theta,\psi))$ .  The curves $\theta=\theta_0$ and $\psi=\psi_0$ for constants $\theta_0$ and $\psi_0$ should produce curves like the black circles in figure \ref{torus}.
\item Find an algebraic equation in $x_1, x_2, x_3$ that defines this same torus of revolution.  ({\em Hint: start with the equation for the circle that gets revolved, using coordinates $y$ and $x_3$.  Then replace $y$ with $\sqrt{x_1^2+x_2^2}$ and get rid of the square root sign.})
\item For a suitably chosen cross section (viewing the picture cut along the $x_2x_3$-plane), the plane $x_1=\frac{r}{\sqrt{R^2-r^2}}x_3$ and the torus described in this problem appear as in figure \ref{hopfcut}.  Can you say anything about the curves which are created by the intersection of this plane and the torus?
\end{enumerate}
\end{exercise}

A torus is often described topologically by identifying the opposite edges of a square.  The old video game {\em asteroids} was played on a torus; when your spaceship left the top of the screen, it reappeared on the bottom, and when it left the left edge, it reappeared on the right.  The torus can be given a geometry (a notion of straight lines, distances and angles) coming from this square.  So one difference between the geometry on a sphere and the geometry on a torus is that on a torus, the sum of the angles of a triangle is $\pi$, while on a sphere, it is always greater.  

We can't preserve this geometry when putting the torus into $\R^3$, but we can if we put it into $\R^4$.  In $\R^3$ we can roll the square into a cylinder, keeping the geometry intact, but then we can't join the ends of the cylinder without stretching the square in places.  With an extra dimension to work with, this problem goes away.

\subsection{$S^3$: the three sphere}

Visualizing the three dimensional sphere can be tricky, but there are many ways to do it, and with practice, you can become quite familiar with this object.

Just as stereographic projection relates the circle to the line, and $S^2$ to the Euclidean plane, stereographic projection gives us a map of all but one point of $S^3$ via a distorted metric on a three dimensional Euclidean space.  One way to get to know $S^3$ is to get to know this distortion.

Stereographic projection involves a choice of a point to project from, and an antipodal point where we imagine the image space touching the sphere.  A pair of antipodal points defines an equator, the sphere of one dimension less that sits halfway between them.  For $S^2$, this is the familiar equatorial circle, and it is sent to the unit circle in the plane via the stereographic projection defined above.  For $S^3$, the equator is a $S^2$, and under a stereographic projection to $\R^3$, it is sent to the unit sphere.

Along this equatorial sphere, stereographic projection does not distort.  Distances and shapes living on this sphere appear in the projection as they are on the sphere.  Outside of the equator, things are stretched.  In the case of $\sigma:S^3\setminus\{\infty\}\to\R^3$, the whole infinite extent of $\R^3$ outside the unit sphere is used to represent only one hemisphere of $S^3$.  In every direction, radiating away from the equatorial sphere, there are lines that seem to diverge, but interpreted in $S^3$, these lines all converge on the north pole. Inside the equatorial $S^2$ on the other hand, everything appears compressed.  Moving with speeds relevant to $S^3$, and not $\R^3$, it takes as long to cross through the interior of the equator as it does to cross the whole of space outside this ball, pass through infinity, and return from the other direction.

The formula for stereographic projection that sends the sphere 
$$
x_1^2+x_2^2+x_3^2+(x_4-1/2)^2=1/4
$$
to the plane $x_4=0$ by projection from $(0,0,0,1)$ is very similar to the formula in one dimension less.  No confusion should arise from calling them both by the same name.
$$
\sigma(x_1,x_2,x_3,x_4)=\left(\frac{x_1}{1-x_4},\frac{x_2}{1-x_4},\frac{x_3}{1-x_4}\right)
$$
Observe that this same formula works to project from the unit sphere, just as the formula in $\R^3$ worked for either the sphere resting on the plane or cut along the equator.

The inverse map sending the plane back to the unit sphere
\\
\\
$\sigma^{-1}(y_1,y_2,y_3)=$ 
$$
\left(\frac{2y_1}{1+\sum_{i=1}^3y_i^2},\frac{2y_2}{1+\sum_{i=1}^3y_i^2},\frac{2y_3}{1+\sum_{i=1}^3y_i^2}, \frac{-1+\sum_{i=1}^3y_i^2}{1+\sum_{i=1}^3y_i^2}\right).
$$

\begin{exercise} Use the fact that the stereographic projection of a circle is a circle or a line to show that in every dimension, spheres are sent to spheres and linear spaces by stereographic projection.
\end{exercise}

\begin{figure}
\centerline{\resizebox{\textwidth}{!}{\includegraphics{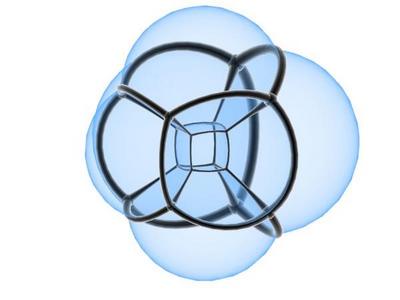}}}
\caption{ A hypercube }
\label{hypercube}
\end{figure}

\begin{exercise}To help us visualize $S^3$, we'll draw some shapes and see what they look like in $\R^3$.  
\begin{enumerate}
\item The vertices of a square inscribed in the unit circle with sides parallel to the axes are the four points $\left\{\left(\frac{\pm1}{\sqrt2},\frac{\pm1}{\sqrt2}\right)\right\}$.  The vertices of a cube similarly inscribed in the unit sphere in $\R^3$ are the eight points  $\left\{\left(\frac{\pm1}{\sqrt3},\frac{\pm1}{\sqrt3},\frac{\pm1}{\sqrt3}\right)\right\}$.  By analogy, the vertices of a hypercube inscribed in the unit $S^3$ in $\R^4$ are the sixteen points $\left\{\left(\frac{\pm1}{2},\frac{\pm1}{2},\frac{\pm1}{2},\frac{\pm1}{2}\right)\right\}$.  Figure \ref{hypercube} shows what the hypercube looks like under stereographic projection.  Describe what is being shown in this picture.  Just like a square has sides that are line segments, and a cube has side that are squared, a hypercube has sides that are cubes.  How many hypercube-side cubes are there in this picture, and how are they configured?

\item The great circles on $S^2$ are useful in understanding the geometry of the two sphere.  What are two possible analogs for the three sphere, and how can they be recognized in the stereographic image?  Give a description similar to the one you gave in exercise \ref{describegreatcircles}.  How could you define the distance between two points on a sphere?
\end{enumerate}
\end{exercise}

\subsubsection{$S^3\subset\C^2$} Since $\R^4$ can also be thought of as $\C^2$, setting 
$$
z_1=x_1+ix_2, \quad z_2=x_3+ix_4,
$$
the equation $x_1^2+x_2^2+x_3^2+x_4^2=1$, rewritten in terms of these complex coordinates, provides the three sphere with the description
$$
S^3=\{(z_1,z_2)\in\C^2\,|\,|z_1|^2+|z_2|^2=1\}.
$$

\begin{exercise}\label{cxline-circle}
\begin{enumerate}
\item Thinking of $S^3$ in this way, consider the surface where $|z_1|^2=|z_2|^2$.   What does this surface look like?  The surface divides $S^3$ into two pieces.  Describe them.
\item What does the intersection of the solutions of the equation $az_1+bz_2=0$ with $S^3$ look like?  (The coefficients $a$ and $b$ can be real or complex.  The surface defined by $az_1+bz_2=0$ is generally called a {\em complex line} because if everything was real instead of complex, it would define a line in the plane.)
\end{enumerate}
\end{exercise}

\subsection{$S^n$}

By now it should be apparent that there are spheres of every dimension.  There is even a reasonable body of mathematics that studies the infinite dimensional sphere, $S^\infty$.  Graphical representations are much harder to come by, since our intuition tends to break down, and it is only with much practice that mathematicians are able to ``see'' in dimensions higher than three or four.  One of the most popular of the higher dimensional spheres is $S^7$, for various reasons that we won't go into in this course.

\newpage

%% file: quaternions.tex
\section{Quaternions}

Points in a one dimensional Euclidan space form a field, the real numbers, and points in a two dimensional space do as well, via the complex numbers.  What about three dimensional space, or higher dimensions?  Are there other fields that naturally correspond to these geometric objects?  It turns out that the answer is no if we insist on looking for fields, but if we relax just a little bit, there are some more rather interesting algebraic structures that can be found.

$\C$ is often defined using $\R^2$ by naming the basis vectors 1 and $i$, and imposing the rule that $i^2=-1$.  We defined the algebraic structure of $\C$ using vector addition and multiplication by similar triangles.  In that development, the fact that $-1$ has a square root takes the form of an observation or a theorem, but it works equivalently well as a definition.  For $\R^4$, it is much harder to draw pictures and make observations based on elementary geometry, so to define an algebraic structure in four dimensions, we'll just give names to the basis vectors and show what they do, and then observe that there is a geometric interpretation.

\begin{defn}
The {\em quaternions} ($\mathbb H$ for short after Hamilton, their discoverer) are the elements of the vector space $\R^4$ with basis $\{1, i, j, k\}$.  Addition is defined by vector addition, and multiplication follows the rules for scalar multiplication by real numbers, and:
$$
i^2=j^2=k^2=-1,
$$
$$
ij=k, \quad jk=i, \quad ki=j,
$$
and
$$
ji=-k, \quad kj=-i, \quad ik=-j.
$$ 
\end{defn}

Because of the last two lines above, it is clear that the quaternions are not commutative with respect to multiplication.  However, the rest of the field axioms do hold for the quaternions, so they are called a {\em skew field}.

\begin{exercise}
Try some quaternionic arithmetic:
\begin{enumerate}
\item $(3-i+2j)+(i-j+2k)$
\item $(2+3i)(1+j-3k)$
\item $(1+j-3k)(2+3i)$
\item $((i+j)2k)(6-j)$
\item $(i+j)(2k(6-j))$
\end{enumerate}
\end{exercise}

Recall that a complex number $z=x+iy$ has a conjugate $\bar z=x-iy$.  Similarly, a quaternion $q=a+bi+cj+dk$ has a conjugate $\bar q=a-bi-cj-dk$.

\begin{exercise}
The complex conjugate is useful when computing the norm and the inverse of a complex number.   Similarly for the quaternionic conjugate.
\begin{enumerate}
\item For a quaternion $q=a+bi+cj+dk$, check that the formula $\sqrt{q\bar q}$ agrees with the Euclidean vector norm $|q|=\sqrt{a^2+b^2+c^2+d^2}$.
\item Inverses can be computed for quaternions using $q^{-1}=\bar q/|q|^2$.   Find $(3+2i-j-k)^{-1}$ using this formula, and check that this is in fact the inverse.
\item It can be checked by hand that $|q_1q_2|=|q_1||q_2|$ but this is a bit tedious.  You can either check this, or compute a few examples so that it seems likely, or just believe it.
\end{enumerate}
\end{exercise}

\subsection{Rotations of $S^2$}

The most important application of the quaternions is to the algebra of rotations in $\R^3$.  It was this application that led Hamilton to discover them.  He was searching for a way to convenient description of rotations allowing for them to be combined, manipulated and performed in sequence.  Standard methods, involving matrices and linear algebra can be cumbersome and opaque.  Since a rotation is specified by an angle and an axis, with the axis is specified by a point on $S^2$, and the angle a point on $S^1$, a rotation is essentially given by three numbers.  The question becomes, given two sets of three numbers determining a pair of rotations, how to write down the triple describing the rotation caused by the composition: first one rotation and then the other.  This problem was extremely vexing to Hamilton and other mathematicians of his time.  One day, so the story goes, Hamilton was walking through Dublin with his wife and insight hit him while crossing the Brougham Bridge.  In order to multiply triples, one needs to actually multiply quadruples.  His description of his moment of understanding captures the electric feeling of instant awareness that is the addictive rush of doing mathematics.
\begin{quote}
{\em That is to say, I then and there felt the galvanic circuit of thought close; and the sparks which fell from it were the fundamental equations between i,j,k; exactly such as I have used them ever since.}
\end{quote}
On being struck by this idea, he carved the equations defining the multiplication for $i$, $j$, and $k$ into the soft stone of the bridge.  The actual inscription is gone, but a plaque has been installed to note this mathematical event.

A complex number can be split into its real and imaginary parts.  Similarly, a quaternion $q=a+bi+cj+dk$ splits into its {\em real} and {\em purely quaternionic} (or simply {\em pure}) parts.  The real part is $a$ and the pure part is $bi+cj+dk$.  A quaternion whose real part is zero is called a {\em pure quaternion}.  The pure quaternions form a three dimensional space, and it is this fact that allows quaternions to be used to describe rotations in $\R^3$.  If $p=xi+yj+zk$ is a pure quaternion corresponding to a point $(x,y,z)\in\R^3$, and $q$ is any nonzero quaternion, compute the conjugate $q^{-1}pq$.  We will show that this defines a rotation of $\R^3$.

\begin{exercise}
For a nonzero quaternion $q$ and a pure quaternion $p$ as above, define the map 
$$
R_q(p)=q^{-1}pq
$$
\begin{enumerate}
\item Show that $R_q(p)$ is also a pure quaternion.
\item Show that $R_q(p)$ is linear: $R_q(\lambda p)=\lambda R_q(p)$ for any real number $\lambda$ and $R_q(p+p')=R_q(p)+R_q(p')$ for any two pure quaternions $p$ and $p'$.
\item Note that $R_q(p)$ is the same rotation as $R_{\lambda q}(p)$ for any nonzero real number $\lambda$.  So we can restrict ourselves to considering rotations defined by the {\em unit quaternions}, those with $|q|=1$.  
\item Show that $|R_q(p)|=|p|$.
\item Show that if $r=bi+cj+dk$ is the purely quaternionic part of $q$, then $R_q(r)=r$.  (It's the axis of rotation!)
\item The plane perpendicular to $r$ is defined by the equation $bx+cy+dz=0$.  Use the fact that $R_q(p)$ is linear and $R_q(r)=r$ to show that this plane is preserved by $R_q$.
\item Choose a $p$ on this plane perpendicular to the axis.  For example, if $b$ and $c$ are not both zero, you can use $p=ci-bj$.  The angle between $p$ and $R_q(p)$ can be computed using the formula 
$$
\cos\theta=\frac{p\cdot R_q(p)}{|p||R_q(p)|}
$$
(The dot in the numerator is the scalar product (a.k.a. the inner product or dot product) for vectors.)  Show that the right hand side is equal to $a^2-b^2-c^2-d^2$, or $2a^2-1$ if $|q|=1$.  This gives the angle of rotation, up to sign.  
\item To find the sign of the rotation, we need to find out if $\{r,p,R_q(p)\}$ is a right handed set of vectors or a left handed set.  Check that if $q=\frac{1+i}{\sqrt2}$ (and hence $r=\frac{i}{\sqrt{2}}$) and $p=j$, that the triple $\{r,p,R_q(p)\}$ is a right handed triple.  This is the case in general.
\item Observe that what we have shown is that for any quaternion, we can write $q=s(\cos\theta+u\sin\theta)$, where $s\in\R$, $u$ is a unit vector in the $\R^3$ of pure quaternions, and the rotation $R_q$ is the rotation about $u$ by the angle $2\theta$.
\end{enumerate}
\end{exercise}

This is an extremely useful way to compute rotations, and it is frequently used in graphics programming for games and computer animations.  Of course, these transformations could also be done using $3\times3$ matrices, but the angle and axis is much easier to read off the quaternion than to extract from the matrix.

\begin{exercise} Show that if $u$ is a pure unit quaternion then $u^2=-1$.  \end{exercise}

\begin{figure}
  \centerline{\resizebox{0.5\textwidth}{!}{\includegraphics{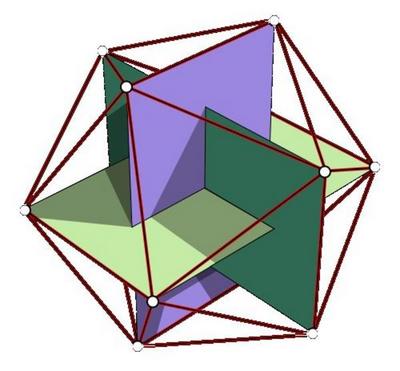}}}

  \caption{The vertices of an icosahedron lie at the vertices of three golden rectangles. ({\em image: Wikimedia})}
\label{icosa}
\end{figure} 

\begin{exercise}\begin{enumerate}
\item Start with an octahedron with vertices on the coordinate axes in $\R^3$.  Rotate the figure about its center so that the edge between $(1,0,0)$ and $(0,1,0)$ becomes vertical.  What $q$ do you use for the rotation, and what are the new coordinates of the vertices?
\item The symmetry group of the octahedron is generated by a rotation of $180^\circ$ around an axis passing through the midpoint of an edge, and a rotation of $120^\circ$ around an axis passing through the center of a face.  Find the quaternions for these generators.
\item The vertices of an icosahedron sit on three golden rectangles lying in orthogonal planes as in figure \ref{icosa}.  (A golden rectangle is one whose sides are in the golden ratio, $1:\frac{1+\sqrt5}{2}$.)  The symmetries of the icosahedron are similarly generated by a $180^\circ$ rotation about an edge midpoint and a $120^\circ$ rotation about a face center.  Find the quaternion generators for this group.
\end{enumerate}
\end{exercise}
 
Perhaps the most exciting part about the description of rotations using quaternions is that the composition of rotations corresponds to quaternion multiplication.

\begin{exercise}
If $q$ and $q'$ are nonzero quaternions, show that 
$$
R_{q'}(R_q(p))=R_{qq'}(p).
$$
\end{exercise}

Observe that since the unit quaternions satisfy the equation $a^2+b^2+c^2+d^2=1$, they form a three dimensional sphere.  Also, if $q$ and $q'$ are unit quaternions then $qq'$ is also a unit quaternion, so the unit quaternions show us that the three dimensional sphere has the structure of a group!  Because we can use any quaternion to give us a rotation via $R_q$, and composition of rotations corresponds to quaternion multiplication, this group is almost $SO(3)$, the group of rotations in $\R^3$.  ($SO$ stands for {\em special orthogonal}, meaning that a rotation is a linear transformation that takes an orthogonal set of vectors or lines, like the coordinate axes, to another set of orthogonal vectors, and that among those, the rotations are special in that they don't change the scale of things or produce any reflections.  So $SO(n)$ stands for the rotations of $\R^n$, or the isometries of the $(n-1)$-sphere.) 

\begin{exercise}
Show that if $q=\sin(\theta/2)+u\cos(\theta/2)$ and $q'=\sin(\psi/2)+v\cos(\psi/2)$ with $u\neq v$ unit pure quaternions, and $R_q=R_{q'}$, then $u=-v$ and $\theta=-\psi$, in other words, $q'=-q$.
\end{exercise}

So every rotation can be represented by a unit quaternion, quaternion multiplication corresponds to composition of rotations, and the only way to write a rotation in two different ways is to multiply be $-1$. That is, we get an isomorphism of groups 
$$
S^3/\{1,-1\}\to SO(3).
$$

\subsection{Rotations of $S^3$}\label{rotS3}

The group structure of $S^3$ also helps to describe the rotations of $S^3$ itself.  Both left and right multiplication are isometries of $S^3$; for each unit quaternion $g,h\in S^3$, the maps
$$
\lambda_g(q)=gq, \quad \rho_h(q)=qh,
$$
denoting left and right multiplication respectively, preserve angles and thus spherical distance.

\begin{exercise}
Check that the angle between $q=\frac{1-3i+k}{\sqrt{11}}$ and $q'=\frac{i-j}{\sqrt2}$ is the same as the angle between $\lambda_k(q)$ and $\lambda_k(q')$.  This is an example of the general fact just stated.
\end{exercise}

Since quaternionic multiplication is not commutative, these two actions are generally distinct, and in fact, the only time when there exist a $g$ and an $h$ so that $\lambda_g(q)=\rho_h(q)$ for all $q$ is when $g=h=\pm1$.  So there is a homomorphism from two copies of $S^3$ to $SO(4)$, the group of rotations of $\R^4$ and the isometries of $S^3$.  The kernel of this map is the group of order two where $g=h=\pm1$, and in fact, this gives all the rotations of $\R^4$; there is an isomorphism of groups
$$
S^3\times S^3/\{(1,1),(-1,-1)\}\to SO(4).
$$

\subsection{Quaternions and $\C^2$}

Another way to think of $\R^4$ and thus the quaternions is as $\C^2$.  If $\R^4$ has coordinates $\{x_1,x_2,x_3,x_4\}$, set $z_1=x_1+ix_2$ and $z_2=x_3+ix_4$.  We could try to ``complexify'' $\C^2$, just like $\R^2$ was complexified (by calling one basis vector 1 and the other $i$ and imposing the rule that $i^2=-1$ on top of scalar multiplication by real numbers).  To construct the quaternions, we can call one basis vector of $\C^2$ 1 and the other $j$, impose the rule $j^2=-1$ on top of scalar multiplication by {\em complex} numbers, and call $ij=k$.  If we also want to have $k^2=-1$, we see that $i$ and $j$ cannot commute, since if they did, then we would have $(ij)^2=ijij=iijj=(-1)(-1)=1$, exactly the opposite of what we want!  So we are led to impose the rule that $ij=-ji$.

Alternatively, we could rethink how we defined the complex numbers, starting with $\R$.  The geometric procedure we used to define multiplication was essentially to define a complex number to be a linear transformation of $\R^2$ that (except for the complex number 0) is bijective and orientation preserving, i.e. one with positive determinant.  A linear transformation $T:\R^2\to\R^2$ is identified with the complex number $T(1,0)$.  In other words, we identify the complex number $z=a+ib$ with the matrix $\twobytwo{a}{-b}{b}{a}$.  Note that the determinant of this matrix is $|z|$.  By analogy, we can consider a linear transformation of $\C^2$ given by the matrix $\twobytwo{z_1}{\bar z_2}{-z_2}{\bar z_1}$ as a representative of the quaternion $q=x_1+ix_2+jx_3+kx_4$.
\begin{exercise}
Check that the matrix sum and product agrees with the quaternion sum and product.
\end{exercise}

Just as the determinant of the matrix $\twobytwo{a}{-b}{b}{a}$ gives $|z|$, the determinant of $\twobytwo{z_1}{\bar z_2}{-z_2}{\bar z_1}$ gives $|q|$.  Since determinants are multiplicative, this shows that $|q_1q_2|=|q_1||q_2|$. 

The complex numbers of unit norm correspond to rotations of $\R^2$.  By analogy, the linear transformations given by matrices $\twobytwo{z_1}{\bar z_2}{-z_2}{\bar z_1}$ with determinant 1 can be thought of as ``rotations'' of $\C^2$.  The name for this group of linear transformations is $SU(2)$, $SU$ for {\em special unitary}.  The fact that $SU(2)$ is almost the same as $SO(3)$ is rather amazing, because in general, these structures are quite different.

\subsection{Octonions}

You might wonder if we could make a similar construction with ${\mathbb H}^2$.  It turns out that you can, and you get an eight dimensional structure called the {\em octonions} or the {\em Cayley numbers} (although Hamilton's friend John Graves discovered them before Cayley).  Like quaternion multiplication, octonion multiplication is not commutative, and on top of that, it is not even associative!  They do retain the rest of the field axioms though.  In particular, each octonion has a multiplicative inverse.  Beyond the octonions, there is no way to extend this process and end up with an algebraic structure where everything has an inverse.  The construction of the octonions and the other details of this discussion would take us too far afield, so you'll have to look it up for yourselves.  There is a very nice book by John Conway and Derek Smith on quaternions and octonions \cite{ConSmi} that you could look at, or there are the articles on John Baez's website \cite{baez}.

\newpage

%% file: hopf_fibration.tex
\section{The Hopf fibration}

So now we are coming to the climax of this wandering investigation.  As suggested at the beginning, not everything that we have seen so far will be directly relevant to understanding the Hopf fibration, but the familiarity that you have developed with $\C$, $\mathbb H$, $S^3$, and so on will give you the ability to see this structure more easily.  Interesting simply for it's esthetic merits, the Hopf fibration rests well balanced at a maximal point of complexity that can be clearly depicted in a drawing or sculpture.  It is also favorite example of many mathematicians who have little interest in making a visual model such as the paintings by Lun-Yi Tsai, as it serves as a launching point into and an illustration of the connections between many fields, ranging from Lie groups and mathematical physics to algebraic topology and homotopy theory.

\begin{figure}
\centerline{\resizebox{\textwidth}{!}{ \includegraphics{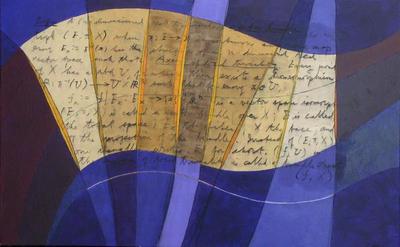}}}
\caption{{\em Purple Vector Bundle} by Lun-Yi Tsai\newline This painting represents a fibration.}
\end{figure}

The Hopf fibration is a mapping from the three dimensional sphere $S^3$ to the two sphere $S^2$.  It presents a window to a deeper understanding of both of these fundamental objects.  A fibration is a very special type of mapping.  Fibrations combine two spaces into a third.  These two spaces are called the {\em base} and the {\em fiber}, and the combination is called the {\em total space of the fibration}, or just the {\em total space} or the {\em fibration}.  Given two spaces $X$ and $Y$, a {\em fibration with base $X$, fiber $Y$, and total space $Z$} is defined using an {\em atlas} $\{U_\alpha\}$ for $X$.  (An atlas for a space $X$ is basically a collection of open sets that completely covers the space, like the paper maps in a world atlas on your bookshelf completely cover the earth.)  Essentially, the total space is defined by giving an atlas for $Z$ where each chart looks like $U_\alpha\times Y$ in a way that is consistent when passing from one chart to the next.  The mapping that is called the fibration is the mapping from $Z$ to $X$ that takes a point represented on a chart $U_\alpha\times Y$ by a pair $(x,y)$ to the point $x$.  For each $x\in X$, there is a copy of $Y$ in $Z$ given by $\{x\}\times Y$.  This is called the {\em fiber over $x$}.  For any pair of spaces, we can define the {\em trivial fibration} where $Z=X\times Y$, and we only need one open set in our atlas to describe the fibration.

\begin{ex}  The M\"obius strip is the simplest example of a non-trivial fibration.  Define the M\"obius strip $M$ as the rectangle $[0,2\pi]\times[-1,1]$ with the relation that $(0,y)=(2\pi,-y)$.  To show that $M$ is a fibration with base $S^1$ and fiber the interval $[-1,1]$, we can use two charts on $S^1$: $U_1=(0,\frac{3\pi}{2})$ and $U_2=(\pi,\pi/2)$.  Note that $U_2$ includes the point $0=2\pi$.  The intersection $U_1\cap U_2$ consists  of  two components, $A=(0,pi/2)$ and $B=(\pi,\frac{3\pi}{2})$.  On $A$, we can use the identification that $(\theta, t)\in U_1\times[-1,1]=(\psi,s)\in U_2\times [-1,1]$ and on $B$, the identification $(\theta, t)\in U_1\times[-1,1]=(\psi,-s)\in U_2\times [-1,1]$.
\end{ex}

\subsection{Hopf fibration via the Riemann sphere}

To begin the description of the Hopf fibration, we look again at $S^3$, described as sitting in $\C^2$:
$$
S^3=\{(z_1,z_2)\in\C^2\,|\,|z_1|^2+|z_2|^2=1\}.
$$
Now, consider the ratio $z_2/z_1$.  This is a complex number, unless $z_1=0$.  But in this case, since we know about $\Pp^1_\C$, the Riemann sphere, by setting $z_2/0=\infty$, we get a mapping
$$
f:S^3\to S^2, \qquad (z_1,z_2)\mapsto z_2/z_1.
$$
This is the Hopf fibration.

\begin{exercise}
Write the Hopf fibration in coordinates as a mapping from $S^3\subset\R^4$ to $S^2\subset\R^3$.  
\end{exercise}

We would like to understand the fibers of this mapping.  For this, it helps to represent $z_1$ and $z_2$ in polar coordinates, so write $z_1=r_1e^{i\theta_1}$ and $z_2=r_2e^{i\theta_2}$. Now observe that for any complex number $\lambda\in\C$ of unit norm ($|\lambda|=1$), and any point $(z_1,z_2)$ in $S^3$, not only is $(\lambda z_1, \lambda z_2)$ still in $S^3$ (since $|\lambda z|=|\lambda||z|=|z|$ for any complex number $z$) but in fact, these points are all on the same fiber of the Hopf fibration; we have $f(z_1,z_2)=f(\lambda z_1, \lambda z_2)$, as
$$\begin{array}{rcl}
f(\lambda z_1, \lambda z_2) & = & \frac{\lambda z_2}{\lambda z_1}\\
&=& \frac{z_2}{z_1}\\
&=& f(z_1,z_2).
\end{array}
$$
The next exercise shows that this is the only way for two points to be on the same fiber.  

\begin{exercise} Suppose that $f(z_1,z_2)=f(w_1,w_2)$  As above, write these points in polar coordinates: $z_1=r_1e^{i\theta_1}$ and $z_2=r_2e^{i\theta_2}$,
$w_1=s_1e^{i\fee_1}$ and $w_2=s_2e^{i\fee_2}$.  

Write $\frac{z_2}{z_1}=\frac{w_2}{w_1}$ as two equations using polar coordinates.  Combine these with the equations coming from the fact that $(z_1,z_2)$ and $(w_1,w_2)$ are on the sphere to show that there exists a unit complex number $\lambda$ such that
$$
(w_1,w_2)=(\lambda z_1, \lambda z_2).
$$
\end{exercise}

Since $(z_1, z_2)$ and $(\lambda z_1, \lambda z_2)$ are distinct points for $\lambda\neq1$, the fibers are topological circles.  They are in fact geometric circles, the same circles that were discovered in exercise \ref{cxline-circle}!

\begin{exercise}
Show that the fibers of the Hopf fibration are the great circles found in exercise \ref{cxline-circle}.

\end{exercise}
\begin{figure}
  \centerline{\resizebox{0.8\textwidth}{!}{\includegraphics{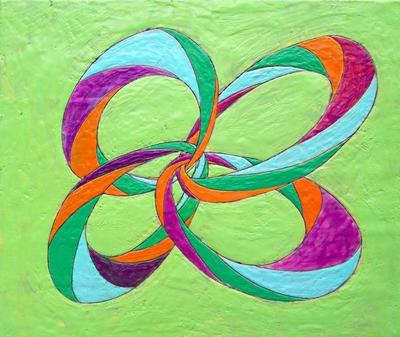}}}
	\caption{{\em Transition Gadget No. 7} \newline Lun-Yi Tsai}
\label{gadget7}
\end{figure}

\subsection{Hopf fibration via quaternions}

Another way of creating the Hopf fibration is by using the $S^3$ of unit quaternions to rotate $S^2$.  If we choose a point $p\in S^2$, then for any quaternion $q$, $R_q(p)$ is also in $S^2$.  So we can define a map $g_p(q)=R_q(p)$, taking $S^3$ to $S^2$.  That is, the image of the point $q\in S^3$ is the point on $S^2$ where $p$ is taken by the rotation $R_q$.  

\begin{exercise}
Show that if $p=(1,0,0)$ then $g_p$ is the same as the map $f$ defined using the Riemann sphere!
\end{exercise}

This means that there is not just one Hopf fibration, but there are in fact infinitely many.  There is a Hopf fibration corresponding to each point $p\in S^2$.

\subsection{Fiber circles are linked}

So there are a two sphere's worth of disjoint circles that fit together to make a three sphere.  We certainly don't expect this to be a trivial fibration, so we'd like to know how the circles fit together.

\begin{exercise}
To do this, consider the equatorial $S^2$, call it $E$, of the unit $S^3\subset\R^4$, where $x_4=0$.  
\begin{enumerate}
\item First, observe that the circle along the equator of $E$, where $x_3=0$ is a fiber of the Hopf mapping.  Call this equator of the equator $F$.
\item Now show that every other point of $S^3$ can be connected to a pair of antipodal points on $E$ by some other fiber of the Hopf mapping.  For $(z_1,z_2)=(x_1,x_2,x_3,x_4)\in S^3$, write $z_1=r_1e^{i\theta_1}$ and $z_2=r_2e^{i\theta_2}$.  For $r_2\neq0$, multiply by some unit complex number $\lambda$ to get a point on $E$.  Since this fiber is a great circle, note that the antipodal point on $E$ is in this fiber as well.
\end{enumerate}
\end{exercise}  

Since every circular fiber contains a pair of antipodal points on the equatorial $S^2$, and joins those antipodal points with a section of the circle in the southern hemisphere of $S^3$ (Under stereographic projection, these points are inside the equatorial $S^2$, where $x_4<0$) and a section in the northern hemisphere (outside $E$ under stereographic projection, where $x_4>0$) every circle in the fibration is linked with $F$.  Just as a two dimensional sphere can be rotated so that any great circle is horizontal, and can be used as an ``equator'', or any pair of antipodal points can be used as the poles, the three sphere can be rotated so that any circle sits where this circle $F$ was sitting.  So actually, any two circles that are orbits of the Hopf fibration are linked, since one can be moved to sit at $x_3=x_4=0$, and then all the rest of the orbits are linked to this particular circle.

\subsection{Latitudinal tori}

\begin{figure}
  \centerline{\resizebox{0.8\textwidth}{!}{\includegraphics{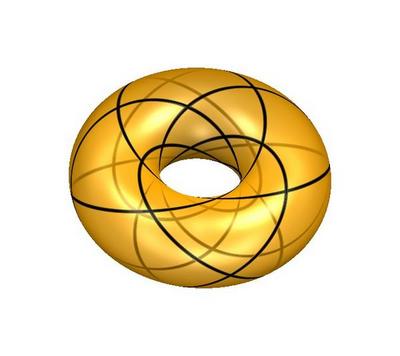}}}
	\caption{Some circles on the torus.}
\label{toruscircles}
\end{figure}

\begin{figure}
  \centerline{\resizebox{0.8\textwidth}{!}{\includegraphics{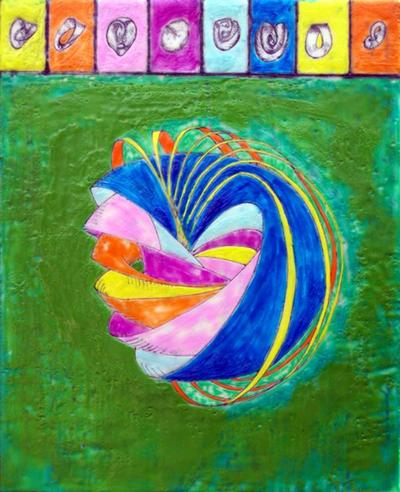}}}
	\caption{{\em Transition Gadget No. 12} \newline Lun-Yi Tsai}
\label{gadget12}
\end{figure}

The fiber circles fit together in another interesting way.  For $z=\frac{z_2}{z_1}\in\Pp^1_\C$, consider the circles $\{|z|=r\}$ for $0<r<\infty$, and the points $z=0$ and $z=\infty$.  These are the lines of latitude on the Riemann sphere, and the poles.  Each of the poles has as its corresponding orbit a circle in $S^3$, and each line of latitude has a circle's worth of orbit circles.
\begin{exercise} Note that these are the tori of exercise \ref{cxline-circle}.  We will use stereographic projection to see how they fit together.
\begin{enumerate}
\item What is $\sigma(f^{-1}(0)$?  How about $\sigma(f^{-1}(\infty)$?
\item You found in exercise \ref{toruseqns} that the equation for a torus of revolution is 
$$
(x_1^2+x_2^2+x_3^2+R^2-r^2)^2=4R^2(x_1^2+x_2^2)
$$
Use this to show that the tori $\sigma(f^{-1}(|z|=r))$ are tori of revolution.  What are the large and small radii?
\item Use what you have found to show that the latitudinal tori are nested as depicted in the paintings by Lun-Yi Tsai, as the tori of revolution generated by revolving a hyperbolic family of Apollonian circles about the line in the family.
\item Show that the fibers of the Hopf fibration are arranged on these tori as the circles in exercise \ref{toruseqns} or figures \ref{toruscircles} and \ref{gadget12}.
\end{enumerate}
\end{exercise}

\subsection{Left and right handed Hopf fibrations}

\begin{figure}
\centerline{\resizebox{.8\textwidth}{!}{\begin{tabular}{c}\includegraphics{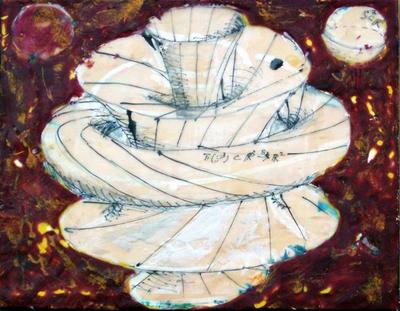} \\ \\ \includegraphics{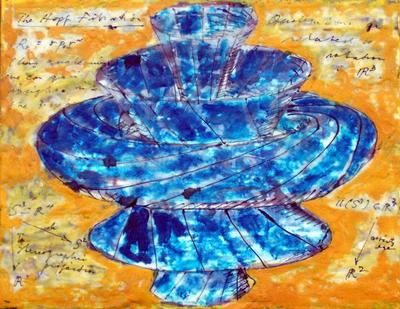}\end{tabular}}}
\caption{ Right and left handed Hopf Fibrations \newline {\em Hopf Fibration Triptych I (top) and II (bottom)} \newline Lun-Yi Tsai}
\end{figure}

In the discussion of quaternions, it was observed that $S^3$ can act on itself in two ways, on the right and on the left.  Related, when we defined the rotation $R_q$, we made a lexicographic choice that had consequences in the geometry of the resulting Hopf fibration because of the non-commuativity of the quaternions.  Instead of $q^{-1}pq$, we could have looked at $qpq^{-1}$.  This is also a rotation of $S^2$.  We can define $L_q(p)=qpq^{-1}$ (note that $L_q=R_{\bar q}$), and look at the Hopf fibrations defined by $h_p(q)=L_q(p)$ for various $p$.

\begin{exercise}
\begin{enumerate}
\item Fix $p=i$ and show that the latitudinal tori for $g_p$ are the same as those for $h_p$.
\item Show that on each of these latitudinal tori, the fibers of $g_p$ wrap around like right handed threads, and the fibers of $h_p$ wrap around like left handed threads.
\end{enumerate}
\end{exercise} 

\subsection{Applications and further directions}

We have only scratched the surface of what can be said about the Hopf fibration.  Hopefully your curiosity has been aroused and you will do some exploring on your own.  Here are some places that you could go to build upon what you already know:
\begin{itemize}
\item Learn about Lie groups and flows by studying the what happens to a $S^3$ when you use the quaternions from a fiber of a Hopf fibration to move it around as in \ref{rotS3}.
\item The study of elementary particles in physics makes use of symmetry groups such as $SU(2)$.  Learn about spin and how this relates to the Hopf fibration.
\item Homotopy groups are the basic objects of study in the field of mathematics called Algebraic Topology.  The fact that the Hopf fibration is a non-trivial fibration $S^3\to S^2$ says something very important about homotopy.  Learn what this is.
\item There are lots of ways that you could make computer generated images of the Hopf fibration using stereographic projection.  Figure out how to make some interesting images or animations.
\item Make a real 3d model of some of the fibers in the Hopf fibration.
\item And so on, there's lots more you could do...
\end{itemize}
